\magnification = 1200
\baselineskip 15pt
\font\tensmc=cmcsc10
\font\bigbf=cmbx10 scaled 1440
\def\mpr#1{\;\smash{\mathop{\hbox to 20pt{\rightarrowfill}}\limits^{#1}}\;}
\def\mpl#1{\;\smash{\mathop{\hbox to 20pt{\leftarrowfill}}\limits^{#1}}\;}
\def\mpd#1{\big\downarrow\rlap{$\vcenter{\hbox{$\scriptstyle#1$}}$}}
\def\mpu#1{\big\uparrow\rlap{$\vcenter{\hbox{$\scriptstyle#1$}}$}}
\font\sf=cmss10
\def\Ct{{\kern1.5pt\hbox{\sf C\kern-5pt I\kern5pt}\kern-2.5pt}}
\def\Qt{{\kern 1.5pt\hbox{\sf Q\kern-6pt I\kern6pt}\kern -2pt}}
\def\Zt{{\hbox{\kern1.5pt\sf Z\kern-8.5pt Z\kern2.5pt}}}
\def\Nt{{\hbox{\sf N\kern-8.5pt I\kern8.5pt}}}
\def\Rt{{\hbox{\sf \kern1pt R\kern-8.2pt I\kern6.8pt}}}
\def\Pt{{\hbox{\sf \kern1pt P\kern-8.2pt I\kern6.8pt}}}
\font\tensmc=cmcsc10

\def\p{\quad$\vbox{\hrule\hbox{\vrule\vbox{\phantom{O}}\vrule}\hrule}$}

\def\s{\vskip 3pt}
\def\cite#1{[#1]}
\def\head#1{\vskip 10pt{\hfil\tensmc #1}\s}
\def\proclaim#1#2{\s{\noindent\bf #1 }{\it #2}\s}
\def\remark#1{\s{\noindent\it #1}}
\def\spa#1{{\tensmc #1}}
\def\title#1{\centerline{\bigbf #1}}

\def\*{\bullet}
\def\N#1{N^{#1}(V_\*)}
\def\n#1{{\bf\underline N}^{#1}(\V)}
\def\V{{\cal U}}
\def\Dt{{\bf D}}
\def\binom#1#2{\left(\matrix{#1\cr #2}\right)}
\def\cl{{\rm cl}}

\vglue 2truecm
\title{ A canonical lift of Chern-Mather classes}
\title{to intersection homology}
\medskip
\centerline{\parindent 10pt J.P. Brasselet\footnote{$^\dagger$}{IML, CNRS,
Luminy Case 930, 132888
Marseille-Cedex 9, France\hfill\break
jpb@iml.univ-mrs.fr } and A. Weber\footnote{$^\ddagger$}{Institute of
Mathematics, Warsaw University, ul. Banacha 2, 02-097 Warszawa, Poland
\hfill\break
aweber@mimuw.edu.pl
\hfill\break\
partially supported by CNRS and KBN grants.\hfill }}
\vskip 40pt

It is well known that for an $n$-dimensional algebraic complex
variety $X$, the Poincar\'e morphism $H^{2n-i}(X) \longrightarrow H_i(X)$,
cap-product by the fundamental class $[X]$, is an isomorphism
if $X$ is a manifold but, in general, it is not. There is
a factorization
$$\matrix{
H^{2n-i}(X) & \mpr{\cap [X]} &H_i(X)\cr
\qquad\alpha\searrow & &\nearrow\omega\qquad \cr
& IH_i(X)&\cr}
$$
by intersection homology groups (we will use only
middle perversity) and intersection homology is the good theory for
considering intersection product of cycles.

On another hand the Chern classes for singular varieties have been defined
by M.H. Schwartz [Sc] and by R. MacPherson [MP1], they are defined in
homology and in general it is not possible to lift them to cohomology.
A natural question arose~: is it possible to lift the
Chern-Schwartz-MacPherson classes to intersection homology?

Firstly J.L. Verdier [V] gave the example of a
singular variety $X$ such that the Chern-Schwartz-MacPherson class
$c_1(X)\in H_2(X)$ can be lifted to $IH_2(X)$ as two distinct Chern
classes of small resolutions $X_j$ of $X$ such that $H_2(X_1) \buildrel
{\cong} \over \rightarrow IH_2(X) \buildrel {\cong} \over
\leftarrow H_2(X_2)$. The computation shows that if we want to express classes
of
singular varieties using small resolution, we need
correction terms living not only in intersection
homology of the singular part.

The second counter-example, due to M. Goresky, is an example of a
singular algebraic variety such that the Chern-Schwartz-MacPherson
class is not in the image of intersection homology, using integer
coefficients. In [BG] was explained the fact that both examples of
Verdier and Goresky are examples of Thom spaces associated to Segre and
Veronese embeddings respectively : ${\Pt}^1 \times {\Pt}^1
\hookrightarrow {\Pt}^3$ and ${\Pt}^2 \hookrightarrow {\Pt}^5$.

The question became : is it possible to lift Chern-Schwartz-MacPherson
classes to intersection homology with rational coefficients ? That time
the answer was positive~: firstly Yokura [Y] proved the result for isolated
singularities, then [BBFGK] proved that all algebraic cycles (and
in particular Chern-Schwartz-MacPherson cycles) can be lifted to intersection
homology, for the middle perversity and with rational coefficients (see
\cite{We} for a simpler proof). Unfortunately, the lifting is not unique, in
general. Due to the Verdier example, it is not obvious that there exists a
canonical lifting.

Among ingredients of MacPherson construction are the Chern-Mather classes, in
fact MacPherson classes are combination of Mather ones. Zhou [Z] proved
existence of lifting of Chern-Mather classes in intersection homology for
total perversity. In this paper, we show, that there exists
a canonical lifting of Mather classes to intersection homology. 
The idea of the proof is the following~: The Chern--Mather classes are
represented by polar varieties. Such polar
variety can be considered as an element of a sequence of
inclusions of polar varieties. The inclusions are of codimension one and
in this case there exists an unique lifting at each step.
As an application,
we obtain canonical lift of Chern-Schwartz-MacPherson classes to intersection
homology for isolated singularities.

\head{1. Recollection of facts about the geometry of Grasmannians}

Let $G(n,m)$ be the Grassmannian of $n$--dimensional spaces in $\Ct^m$ and
let $$V_\*=\{V_0=\{0\}\subset V_1\subset V_2 \subset \dots \subset
V_m=\Ct^m\}$$
be a flag. The flag manifold, set of such flags, will be denoted by ${\bf
F}(m)$.
Define for $i\geq 0$ the Schubert variety \cite{Ch}, \cite{Eh}
$$M^i(V_\*)=\{W\in G(n,m): W+V_{m-n+i-1}\not= \Ct^m\}$$ 
of codimension $i$ in $G(n,m)$.
The cycle $(-1)^iM^i(V_\*)$ represents the Poincar\'e dual of the $i$--th
Chern class $c^i\in H^{2i}(G(m,n))$.
Each $M^i(V_\*)$ has a natural stratification whose smooth strata are:
$$M^{i,k}(V_\*)=\{W\in G(n,m):codim\,( W+V_{m-n+i-1})=k+1\}\,.$$
The regular stratum of $M^i(V_\*)$ is $M^{i,0}(V_\*)$.

\proclaim{Proposition 1.1.}{The Schubert varieties $M^i(V_\*)$ have the
following properties:\hfill\break
1. $M^{i+1}(V_\*)\subset M^i(V_\*)$;\hfill\break
2. $M^{i+1}(V_\*)\cap M^{i,0}(V_\*)\subset M^{i+1,0}(V_\*)$;\hfill\break
3. for $i< n$ the regular part of $M^{i+1}(V_\*)$ is not contained in the
singularities
of $M^i(V_\*)$;\hfill\break
4. $M^0(V_\*)=G(n,m)$, $M^{n+1}(V_\*)=\emptyset$.}

\remark{Proof.} The proof of 1 and 2 is clear, since $V_{m-n+i-1}\subset
V_{m-n+i}$.

For proving 3, suppose
$V_{m-n+i}=V_{m-n+i-1}+lin\{\alpha\}$ where $lin\{\alpha\}$ 
denotes the line generated by $\alpha\in \Ct^m $. This space is not 
equal to $\Ct^m$ and there
 exists $W\in M^{i,0}(V_\*)$ such that $\alpha
\in W$. Then $dim\,(V_{m-n+i}+W)=dim\,(V_{m-n+i-1}+W)=m-1$ and
$W\in M^{i+1,0}(V_\*)$.

The proof of 4 follows by dimension considerations.\p\s

Since $M^{i+1}(V_\*)$ is irreducible, then a generic point of $M^{i+1}(V_\*)$
is in $M^{i,0}(V_\*)$.

\head{2. The Gauss map and polar varieties}

Let $X^n \subset \Ct^m$ be an affine variety. 
Let us denote by $\Sigma_X$ the singular part of $X$ and by $X_{\rm reg} =
X\setminus\Sigma_X$ the regular one. There is a natural map $s:X_{\rm reg} 
\to G(n,m)\times \Ct^m$ defined by $s(x)=(T_x(X_{\rm reg}),x)$. The 
Nash blowup is the closure of the image of $s$. We have the diagram:
$$\matrix{
&\hat X& \hookrightarrow & G(n,m)\times \Ct^m & \mpr{p_1} & G(n,m)\,,\cr
\hfill ^s\nearrow &\mpd\pi && \mpd{p_2}\cr
X_{\rm reg}\;\hookrightarrow& X&\hookrightarrow&\Ct^m}$$
where $\pi=p_{2|\hat X}$ and $g=p_{1|\hat X}:\hat X \rightarrow G(n,m)$
is the Gauss map. For a general flag
$V_\*$ we define a cycle
$$\N i=\pi[g^{-1}M^i(V_\*)]={\rm closure}\{ \,\pi[g^{-1}(M^i(V_\*))\cap
s(X_{\rm reg})]\,\}.$$
We will state the genericity conditions.

\remark{Definition 2.1.} We say that the map $g:\hat X \rightarrow G(n,m)$ is
general if it is transverse to all strata $M^{i,k}(V_\*)$. This means that g
restricted
to $s(X_{\rm reg})$ and to each stratum of the special fiber is
transverse to the strata $M^{i,k}(V_\*)$.

Let us consider the standard flag $$V^0_\*=\{\{0\}\subset
lin\{e_1\}\subset lin\{e_1,e_2\}\subset\dots\subset \Ct^m\}\,.$$ The
group $Gl(m)$ acts on each $G(n,m)$ transitively.
By Kleinman's theorem (\cite{Kl} 2. Theorem) there
exists an open algebraic subset $U$ of $Gl(m)$ such that for any
$a\in U$ the map $a\cdot g$ is transverse to the strata $M^{i,k}(V_\*)$.
For an open subset $\V\in{\bf F}(m)$ define the total polar variety
$$\n i=\{(x,V_\*)\in X\times {\bf F}(m): x\in \N i, \, V_\*\in\V\}\,.$$
It is an algebraic set over $\V$. The projection on $\V$
is not a fibration in general. We fix a sufficiently small $\V\subset
{\bf F}(m)$ such that the projections $\n i\rightarrow\V$ are
fibrations for all $i\geq 
0$ and $\V\subset U\cdot V^0_\*$ for $U$ as in the Kleinman's theorem.

\remark{Definition 2.2.} The flag $V_\*$ is called {\it good} if
$V_\*\in\V$. 

For any good flag $V_\*\in \V$ the cycle $(-1)^{i}\N i$
is called the polar variety \cite{LT}, \cite{Pi}, it represents
the Chern--Mather class $c^M_{n-i}(X)$ \cite{MP1}.
It is the closure of critical
points of the projection $X_{\rm reg}\rightarrow V_{m-n+i-1}^\perp$~.

\head{3. Construction of a lift of Chern classes to intersection homology}

Let us denote by $\Dt _X$ the dual sheaf (Borel-Moore homology sheaf) of an
algebraic complex variety, and $IC_X$ its intersection homology
sheaf.  All coefficients of homology and intersection homology
are rationals.

The dual proof of \cite {BBFGK}
(\S 3.5) shows that for closed embedding of codimension one $W\hookrightarrow
X$,
there exists  a lift $\nu$ of the natural morphism $\iota:
\Dt _W \mpr{} \Dt_X $
providing a commutative diagram 
$$
\matrix{
IC_W &\mpr{\nu} &IC_X\cr
\downarrow & & \downarrow\cr
\Dt _W&\mpr{\iota} &\Dt _X\cr
} \qquad \hbox{ and then } \qquad
\matrix{
IH_*(W) &\mpr{\nu_*} &IH_*(X)\cr
\downarrow & & \downarrow\cr
H_*(W)&\mpr{\iota_*} &H_*(X)\,.\cr
}
$$
It is unique as soon as there exists an unique lift on the
smooth part of W ({\it loc. cit.} p.
167). Here $IC_W$ is the intersection homology sheaf of $W$
considered as a sheaf on $X$ supported by $W$ and $\Dt_W$ is
the Borel--Moore homology sheaf also considered as a sheaf on $X$.

Let $\N\*$ be the sequence of polar varieties associated to a good 
flag:
$$\N n \subset  \dots \subset \N 1 \subset \N0=X^n\,.$$
Then, for each $j$, $\N j$ has codimension one in $\N{j-1}$ and no
component of $\N j$
is contained in the singularities of $\N{j-1}$.
In this situation, by the previous result, there
exist unique sheaf morphisms 
$$IC_{\N n} \mpr{} \dots \mpr{}IC_{\N 1} \mpr{}  IC_{\N 0} =IC_X$$
which are lifts of the natural morphisms  
$$\Dt _{\N n}\mpr{}  \dots \mpr{} \Dt _{\N 1} \mpr{} \Dt _{\N 0} = \Dt _X\,.$$
We obtain the induced diagram of morphisms of intersection homology
and homology groups: 
$$\matrix{
IH_*(\N n) & \mpr{}& \dots &\mpr{}& IH_*(\N 1)& \mpr{}& IH_*(\N 0) 
& = &IH_*(X)\cr
\downarrow & &  & & \downarrow & &\downarrow & &\downarrow\cr
H_*(\N n) & \mpr{} & \dots &\mpr{}& H_*(\N 1)& \mpr{} & H_*(\N 0)
& = & H_*(X)\,.}$$

The fundamental class of the polar variety $[\N i]$ belongs to
$IH_{2(n-i)}(\N i)$.

\remark{Definition 3.1.} We define an element $\tilde c^i(X) =\tilde
c_{n-i}(X)\in
IH_{2(n-i)}(X)$
as the image (under the considered sequence of morphisms) of the fundamental
class of the polar variety $(-1)^i[\N i]\in IH_{2(n-i)}(\N i)$.

\proclaim{Proposition 3.2.}{The class $\tilde c^i(X)$ does not depend on the
choice of the good flag.}

\remark{Proof.}  
We have two sequences of morphisms: 
$$\matrix{
IH_*(\n n) & \rightarrow & \dots & \rightarrow &IH_*(\n 1) & \rightarrow & 
IH_*(\n 0)& = & IH_*(X)\otimes H^*(\V)\cr
\uparrow &&&&\uparrow & & \uparrow & &\uparrow\cr
IH_*(\N n)  & \rightarrow & \dots & \rightarrow & IH_*(\N 1) & \rightarrow & 
IH_*(\N 0)& = & IH_*(X) \,,}$$
where the top row is written for the total polar variety of good
flags and the bottom row is written for a fixed flag. The diagram commutes 
since
$\N i\hookrightarrow \n i$ is normally nonsingular.
The class $(-1)^i[\N i]\in IH_{2n-2i}(\n i)$ does not depend on the choice of
the
good flag $V_\*\in \V$ and is the class of a fiber in the bundle $\n i
\rightarrow \V$. Applying the sequence of maps of intersection homology groups
we obtain an element $(-1)^i[\N i]\in IH_{2n-2i}(X\times \V)$ which is
independent on $V_\*$ . It can be written as
 $\tilde c^i(X)\otimes [pt]\in IH_{2n-2i}(X)\otimes H_0(\V )$. \p

\head{4. Chern class for quasi--projective variety}

Let $X^n\subset \Pt^m$ be a smooth quasi--projective variety. Let $T_{CX}$ be
the tangent bundle of the affine cone $CX\subset \Ct^{n+1}$
over $X$ away of origin. The bundle $T_{CX}$ is induced from a bundle
$\tau \rightarrow X$. The polar varieties of $X$ are defined to be the
projectivization of
the ones of $CX$. They represent the Chern classes of the bundle $\tau$.
To recover the Chern classes of $X$ we use the following formulas:

\item {1.} a formula for bundles (see \cite{MS~: 14.10} for the case
$X^n=\Pt^n$):
$$T_X\oplus\Theta \simeq \tau\otimes\gamma^*\,;$$
where $\Theta$ is the trivial bundle and $\gamma={\cal O}(-1)$ is
the tautological bundle.

\item {2.} a formula for Chern classes of a tensor product: let $E$
be a $k$--dimensional bundle and let $L$ be a line bundle,
$c_i=c^i(E)$ and $a=c^1(L)$ then

$c^*(E\otimes L) = 1 + c_1 + ka+ \hfill (4.1)$

$\phantom{c_*(E\otimes L)=} c_2+\binom{ k-1} 1 ac_1 + \binom k 2 a^2+$

$\phantom{c_*(E\otimes L)=} c_3+\binom{ k-2} 1 ac_2 + \binom{ k-1} 2 a^2 c_1
+\binom k 3 a^3+$

$\phantom{c_*(E\otimes L)=} c_4+\binom{ k-3} 1 ac_3 + \binom{ k-2} 2 a^2 c_2
+\binom {k-1} 3 a^3c_1+\binom k 4 a^4+\dots$\s

If we put $k=n+1$,  $E=\tau$, $L = \gamma^*$ then 
$a\in H^2(X)$ is the class of hyperplane section and we obtain
a formula for the Chern class of $X$.

Suppose $X$ is singular. There is no tangent bundle to $CX$ nor a bundle
$\tau$. Instead we set $c_i=(-1)^i[\N i]\in H_{2n-2i}(X)$ where 
 $\N i\subset X$ is the projectivization of the polar variety of $CX$.
Then the formula (4.1) computes the Chern--Mather class of $X$, [Pi].
By the same formula we define a lift of the Chern--Mather class to intersection
homology of $X$, but now $c_i$ is the lift of $(-1)^i[\N i]$ to 
$IH_{2n-2i}(X)$ constructed in \S 3.

\head{5. A lift of Chern-Mather classes}

The Proposition 3.2 and the formula 4.1 provides us with a method of
defining a canonical lift of Chern-Mather classes:

\proclaim{Theorem 5.1.}{ The Chern-Mather classes of an algebraic complex
variety can be lifted to intersection homology, in a canonical way.}

The Chern-Schwartz-MacPherson class is a combination:
$$c_*(X)=\sum n_\alpha incl_*c^M_*(\overline S_\alpha)\,,$$
where $\{S_\alpha\}$ is the minimal stratification (see \cite{Te}). Thus we
obtain a canonical lift of Chern-Schwartz-MacPherson classes to $$\bigoplus
IH_*(\overline S_\alpha)\,.$$

Suppose $X$ admits only isolated singularities $\{a_i\}$, then the total
Chern-Schwartz-MacPherson class is equal to~:
$$c_*(X) = c_*^M(X) + \sum (1 - Eu_{a_i}) [a_i]\,. $$
It can be lifted to $IH_*(X)$ as soon as we lift the class
$[a_i]$. Few canonical liftings can be defined but they
coincide if $X$ is irreducible.

\proclaim{Theorem 5.2.}{ The Chern-Schwartz-MacPherson
classes of an irreducible algebraic complex
variety which has only isolated singularities can be lifted
to intersection homology, in a canonical way.}

\head{6. The class $\tilde c^1(X)$ and a small resolution}

Let us recall that a {\it small resolution} $\varpi:\tilde X \rightarrow X$ is
a resolution for which there exists a stratification $S_\alpha$ of $X$ such that
for any $x\in S_\alpha$, $\dim \varpi^{-1}(x) < 1/2\, {\rm codim}
(S_\alpha)$. In this case, there is an identification of perverse
sheaves $R \varpi_*\Qt_{\tilde X}\cong IC_X$ and the intersection homology
groups of $X$ are identified with homology groups of $\tilde X$ ([GM], \S 6.2
and [MP2] \S 5).

We will show that:\s

\proclaim{Proposition 6.1.}{For a small resolution we have:\hfill\break
 1. {The lift $\tilde \N 1$ of the cycle $\N 1$ to $\tilde X$
(proper inverse image) represents the class $-\tilde c^1(X)\in H_{2n-2}(\tilde
X)
\simeq IH_{2n-2}( X)$;}\hfill\break
 2. { If $X$ does not have singularities in  codimension one (e.g. if $X$ is
normal), then $\tilde c^1(X)$ coincides with $c^1(\tilde X) \in
H_{2n-2}(\tilde X)$.
}}

\remark{Proof of 1.} We have the following diagram of sheaves over $X$:\s\s

\hfil$\matrix{R \varpi_*\Qt_{\tilde X}&\mpr{\simeq}&IC_X\cr
\mpu{}& &\mpu{\nu\;\rm unique}\cr
R \varpi_*\Qt_{\tilde \N 1}&\hookleftarrow&IC_{\N 1}}$

\hfil$^{\rm \qquad direct\, summand}$

\noindent We may assume that it commutes away of singularities of $\N 1$ and
$X$. Then it commutes on the whole $X$ since $\nu$ is unique. Thus the class
$-[\tilde \N 1]\in H_{2n-2}(\tilde X)$ corresponds to $\tilde c^1(X) \in
IH_{2n-2}(X)$.\p\s

The result is not true for higher classes (see Observation 7.2).

\remark{Proof of 2.} Suppose that $codim\, \Sigma_X\geq 2$. Then $H^{2n-2}(X)
\rightarrow
IH_2(X)$ is surjective and thus $H_{2n-2}(\tilde X)=IH_{2n-2}(X) \rightarrow
H_{2n-2}(X)$ is injective. In $H_{2n-2}(X)$ we have:
$$ \varpi_*c^1(\tilde X) = c^1(X) + \sum n_\alpha[\overline
S_\alpha^{n-1}]\,.$$
Since $X$ has no singularities in codimension one, then all
$S_\alpha^{n-1}=\emptyset$. Thus $ \varpi_*c^1(\tilde X) = c^1(X) $ in
$H_{2n-2}(X)$.
The induced morphism $\varpi_*$ is injective.
The corresponding equality clearly holds in homology of $\tilde X$.\p

In general the first Chern class $\tilde c^1(X)$ differs from the Chern class
of a small resolution. Let us consider two examples.

\remark{Example 6.2.} Let $X\subset \Pt^2$ be given by an equation
 $f(x,y) =xy=0\,.$
It admits a small resolution which is its normalization
$$\tilde X\simeq\Pt^1_1\sqcup\Pt^1_2\mpr\varpi X\,.$$
The Chern class of $\tilde X$, in $H_*(\tilde X)$, is
$[\tilde X]+2[\tilde{pt}_1]+2[\tilde{pt}_2]$, where $\tilde{pt}_1$ is a point
 in $\Pt^1_1$ and $\tilde{pt}_2$
is a point in $\Pt^1_2$.
Let us compute $\tilde c^1(X)$.
The affine cone $CX$ consists of two planes. The projection of it along
a general direction is nonsingular.
Thus $\N 1=\emptyset$.
To compute $\tilde c^1(X)$ we use the formula 4.1:
$$\tilde c^*(X)=1+(-[\N 1]+2a)+\dots\,,$$
where $a=[pt_1]+[pt_2]$ is the class of a hyperplane section.
We obtain $\tilde c^1(X)=2[pt_1]+2[pt_2]$. In this case the class
$\tilde c^1(X)$ coincides with $c^1(\tilde X)$.
\s

\remark{Example 6.3.} Let $X\subset \Pt^2$ be given by the equation
 $f(x,y,z)=x^3+y^2z=0\,.$ It has an singular point denoted by $\{ x_o\}$. 
It admits a small resolution which is also its normalization:
$$\tilde X\simeq\Pt^1\mpr\varpi X\,.$$
The Chern class of $\tilde X$, in $H_*(\tilde X)$, is $[\tilde X]+2[pt]$.
The class $\tilde c^1(X)$ is the homological Chern--Mather class since
$X$ is a topological manifold.
We compute the Chern--Mather class $c_*^M (X)$ from the formula:
$$\varpi_*(c_*(\tilde X))=c_*^M(X)+k\,incl_*c_*^M(\{x_o\})\,,$$
where $k$ is given by the following
expression in terms of local Euler obstruction \cite{MP1}~:
$$\varpi_*(1_{\tilde X})=Eu_X+k\,Eu_{\{x_o\}}\quad{\rm i.e.}\quad k=-1\,.$$
We obtain $\tilde c^*(X)=c^M_{1-*}(X)=[X]+3[pt]$. In this case the difference
between $\tilde c^1(X)$ and $c^1(\tilde X)$ is the class $[pt]$.
\s
Let us give a general expression of the difference between $\tilde c^1(X)$
and the first Chern class of a small resolution.

We remind that $X\subset \Pt^m$. Let $U\subset \Pt^m$,
$U\simeq\Ct^m$  be one of the standard affine charts. 
Let $W$ be a irreducible component of $\varpi^{-1}\Sigma_X$ with
$dim\,W=dim\,X-1=n-1$. Since $\varpi:\tilde X\rightarrow X$ is a small
resolution thus $dim\,\varpi(W)=n-1$. To each such $W$ we will assign
a number. 

\remark{Definition 6.4.}  Assume that $U\cap
\varpi(W)\not=\emptyset$.  We define the  {\it Jacobian multiplicity} of $W$,
denoted by $n_W$, as the order of zeros on $W$ of the Jacobian of the
composition: $$\varpi^{-1}(U\cap X) \mpr{\varpi} U\cap X\mpr p \Ct^n\,,$$
where $p$ is a general projection from $U\simeq\Ct^m$ to $\Ct^n$.

\remark{Remark 6.5.} For a small resolution $\varpi:\tilde
X\rightarrow X$  let us define by $C^\circ\tilde X$ the pull-back
(fibred product):
$$
\matrix{
C^\circ\tilde X &\mpr{C^\circ\varpi} &C^\circ X&=&CX\setminus \{0\}&\subset&\Ct^{m+1}\cr
\downarrow & & \downarrow\cr
\tilde X&\mpr{\varpi} &X&\subset&\Pt^{m}\,.\cr
}$$
 Then, in the Definition 6.4, instead of a local
general projection $p$ we can take a general global projection
$p:\Ct^{m+1}\rightarrow \Ct^{n+1}$ and  use the composition 
$\tilde p=p\circ C^\circ\varpi$:
$$C^\circ\tilde X\mpr{C^\circ\varpi} C^\circ X\mpr p \Ct^{n+1}$$
to compute the multiplicity of $W$.
\s

The Jacobian multiplicity can be expressed by the local Euler obstruction. 
For $X$
locally irreducible, the Jacobian multiplicity of $W$ equals $Eu_X( x) -1$,
where $x$ is a generic point of $\varpi (W)$. If $X$ is not locally
irreducible,
we should take a suitable local component of $X$. For proving this, it is
sufficient to look at the case where $X$ is a curve.

\proclaim{Theorem 6.6.}{ Let $\varpi:\tilde X\rightarrow X$ be a small
resolution.
Then $c^1(\tilde X)=\tilde c^1(X)-\sum n_W\,[W]$,
where the sum runs over the set of irreducible components of
$\varpi^{-1}(\Sigma_X)$
such that  $dim\,W=n-1$.}

\remark{Remark 6.7.} The corresponding relation between Chern-Mather
 classes in homology can be found in \cite{MP1}.

\remark{Proof.} Firstly notice that the polar variety $C^\circ\tilde\N
1$ of $C^\circ\tilde X$ is the
closure of zeros
of the Jacobian of a generic projection $\tilde p:C^\circ\tilde X\rightarrow
\Ct^{n+1}$.
If the map
$\tilde p$ is not generic, then we should take into account the multiplicity
of the zeros of the Jacobian. 
The singularities of $\tilde p$ consist of the
singularities of $p$ and of the singularities of $C^\circ\varpi$. Thus
the components of the polar variety in $C^\circ\tilde X$ come from
$C^\circ\N 1\subset C^\circ X$ or from the singularities of $C^\circ
X$. If $p$ is general then these two sets of components are disjoint.
The components of $\varpi^{-1}\Sigma_X$ should be counted with multiplicities
$n_W$. The Chern class of $\tilde X$ is
$$c^1(\tilde X)=-([\tilde\N 1]+\sum n_W\,[W])+(n+1)\tilde a=
\tilde c^1(X)-\sum n_W\,[W]\,.$$
where $\tilde a$ is the class of a hyperplane section of
$\tilde X$ which is the inverse image of the class of a
hyperplane section of $X$.
\p
\remark{Explanation of the examples.} In the first case the components of
$\varpi^{-1}\Sigma_X$ are two points, but with zero Jacobian multiplicity.
In the second example consider a general local projection from $X\setminus
\{z\not=0\}$ to $\Ct$.
A point $t$ of normalization $\tilde X\simeq\Pt^1$ is sent to $[t^2:t^3:1]$
and then projected
to a point $at^2+bt^3$. Thus the Jacobian multiplicity is one for $t=0$.

\head{7. The example of J. L. Verdier {\rm (see \cite{V} and \cite{BG})}}

Let $B=\Pt^1_x\times\Pt^1_y \hookrightarrow \Pt^3$ be the Segre embedding:
$$([x_0:x_1],[y_0:y_1])\mapsto [x_0y_0:x_0y_1:x_1y_0:x_1y_1]\,.$$
The quadric $B$ is described by the equation:
$$z_0z_3 - z_1z_2=0\,.$$
Denote by $X$ the projective cone over $B$:
$$X=cB\subset\Pt^4$$
defined by the same equation in $\Pt^4$.
Topologically $X$ is the Thom space of the bundle $\gamma_{|B}$; where
$\gamma$
is the tautological bundle over $\Pt^3$.
The variety $X$ admits two small resolutions. To see them consider the bundle
$\gamma_{|B}\rightarrow B=\Pt^1_x\times\Pt^1_y$ as the family of bundles over
$\Pt_y^1$ parameterized by $\Pt_x^1$. For each $x$ the bundle
$\gamma_{|\{x\}\times\Pt_y^1 }$ is equivalent to the tautological bundle over
$\Pt_y^1$. We apply the construction of Thom space for each $x$. We obtain a
smooth space $X_1$ fibered over $\Pt_x^1$ with fiber $c\, \Pt_y^1\simeq
\Pt^2$. The space $X_1$ is a small resolution of $X$; the inverse image of the
singular point is the set of infinity points of the family of the Thom spaces,
i.e. it is $\Pt^1$. The second small resolution $X_2$ is obtained by changing
the
role of
$x$ and $y$.

We have the canonical isomorphisms:
$$H_*(X_1)\simeq IH_*(X)\simeq H_*(X_2)\,.$$
We will calculate the intersection homology groups of $X$ with
rational coefficients. Since it is the Thom
space, therefore
$$\tilde H_*(X) \simeq H_{*-2}(B) \simeq
0\,, 0\,, \Qt\,, 0\,, \Qt^2\,, 0\,, \Qt\,,$$
and
$$IH_*(X) \simeq \left\{\matrix{
H^{6-*}(X)\quad &{\rm for}& \quad *<3\cr
im\,PD\qquad &{\rm for}&\quad *=3\cr
H_*(X)\qquad &{\rm for}&\quad *>3
}\right.\,,$$
where $PD : H^{3}(X) \to H_3(X)$ is the Poincar\'e homomorphism,
cap-product by the fundamental class $[X]$, thus
$$IH_*(X) \simeq\Qt\,, 0\,, \Qt^2\,, 0\,, \Qt^2\,, 0\,, \Qt\,.$$
We will describe the generators (see the figures 1--6):

$IH_2(X)$ is generated by the projective lines: $[\Pt^1_x]=d_1$ and
$[\Pt^1_y]=d_2$.

$IH_4(X)$ is generated by the cones: $c(d_1)=p_1$ and $c(d_2)= p_2$.

\noindent The corresponding generators in $X_1$ and $X_2$ are the proper
inverse images of those in $X$ and will be denoted by the same letter.
The homological Chern class of $X_1$ and $X_2$ were calculated in \cite{BG}
and
they are the following: $$c^*(X_1)=[X_1]+(3p_1+3p_2)+(3d_1+5d_2)+6\{ pt\}\,,$$
$$c^*(X_2)=[X_2]+(3p_1+3p_2)+(5d_1+3d_2)+6\{ pt\}\,.$$
This shows, that the Chern class of $X$ cannot be calculated using small
resolution without correction terms in $H_2(X_i)$.

Now we will calculate $\tilde c^*(X)$ straightforward.
Firstly we find suitable polar varieties. The cone over $X$ in $\Ct^5$ is
described by the equation:
$$f(z)=z_0z_3 - z_1z_2=0\,.$$
The gradient field of $f$ is
$$grad\,f(z)=(z_3,-z_2,-z_1,z_0,0)\,.$$
Fix the flag
$$V_\*=\{\{0\},\,lin\{e_0-e_3\},\,lin\{e_0-e_3,e_1-e_2\},\,lin\{e_0,e_1-e_2,e_
3\},\,
lin\{e_0,e_1,e_2,e_3\},\Ct^5\}\,.$$
Let $CX^o=C(X_{\rm reg})\setminus \{0\}=CX\setminus\{z_4=0\}$. Then

$C\N 1=\cl\{z\in CX^o: grad\,f(z)^\perp\oplus lin\{e_0-e_3\}\not=\Ct^5\}\,,$

$\phantom{C\N 1}
=\cl\{z\in CX^o: grad\,f(z)\perp (e_0-e_3)\}\,,$

$\phantom{C\N 1}
=\cl\{z\in CX^o: z_0-z_3=0\}\,,$
\s

$C\N 2=\cl\{z\in CX^o: grad\,f(z)^\perp\oplus
lin\{e_0-e_3,e_1-e_2\}\not=\Ct^5\}\,,$

$\phantom{C\N 1}
=\cl\{z\in CX^o: grad\,f(z)\perp (e_0-e_3),\;grad\,f(z)\perp (e_1-e_2)\}\,,$

$\phantom{C\N 1}
=\cl\{z\in CX^o: z_0=z_3,\; z_1=z_2\}\,,$
\s

$C\N 3=\cl\{z\in CX^o: grad\,f(z)^\perp\oplus
lin\{e_0,e_1-e_2,e_3\}\not=\Ct^5\}\,,$

$\phantom{C\N 1}
=\cl\{z\in CX^o: grad\,f(z)\perp e_0,\;grad\,f(z)\perp
(e_1-e_2) ,\;grad\,f(z)\perp e_3\}\,,$

$\phantom{C\N 1}
=\cl\{z\in CX^o: z_0=z_3=0,\;z_1=z_2\}=\emptyset\,,$

\noindent see the figures 7--9. The Chern classes of the bundle $\tau$ over
$X_{\rm reg}$ (see \S4)
are represented by the cycles $(-1)^i\N i$.
The cycle $\N 1$ is allowable in $X$. It is the projective cone over the
hyperplane section of $B$. Thus $\tilde c^1(\tau)=-[\N
1]=-(c(d_1)+c(d_2))=-(p_1+p_2)$. To calculate $\tilde c^2(\tau)$ we have 
firstly to find the class $[\N 2]\in IH_2(\N 1)$. The cycle $[\N 2]$ itself 
is not allowable in $\N 1$. To see how it lifts to $IH_2(\N 1)$ let us examine 
the only
singular point of $\N 1$. In the affine chart $\{z_4\not=0\}$ the set $\N 1$
is described by the equations
$$\{z_0z_3-z_1z_2=0,\;z_0=z_3\}=\{z_3^2-z_1z_2=0,\;z_0=z_3\}\,.$$ It is a
singularity of the type $A_1$.

\proclaim{Fact 7.1.} {The surface in $\Ct^3$ with a singularity of type
$A_1:\{z_3^2-z_1z_2=0\}$ is a rational homology manifold.}
\remark{Proof.} We compute the cohomology of the link $L$ of the singular
point.
We use the Gysin sequence of the fibration $p:L\rightarrow L/S^1=K$, where $K$
is a quadric in $\Pt^2$:\vskip 2pt
$$\matrix{H^1(K)&\mpr{p^*}&H^1(L)&\mpr{\int}&H^0(K)&\mpr{2\cdot}&H^2(K)\cr
\| &&&&\| &&\|\cr
0 &&&&\Qt &&\Qt
}$$
We see that $H^1(L)=H^2(L)^*=0$.\p
\s
The space $\N 1$ is the projective cone over $K$.
It is a rational homology manifold.
The polar variety $\N 2$ is the projective cone over two points in $K$.
The group $IH_2(\N 1)=H_2(\N 1)$ is
generated by $[K]$. Thus $[\N 2]$ is a multiple of $[K]$.
To find the multiplier
we intersect $[\N 2]$ and $[K]$ with $[K]$:

$[\N 2]\cdot [K]=2$,

$[K]\cdot [K]=deg\,K =2$.

\noindent We conclude that $[\N 2]=[K]$ in $H_2(\N 1)=IH_2(\N 1)$.

Now $K$ is allowable in $\N 1$ and in $X$; it is the hyperplane section of
$B$.
Thus $[K]=d_1+d_2$. We find that
$$\tilde c^*(\tau)=[X]-(p_1+p_2)+(d_1+d_2)\,.$$
The class $c^1(\gamma^*)$ is represented by $b=[B]\in H^2(X)$
-- the hyperplane section
of $X$. To calculate the class $\tilde c^*(X)$ we use the formula 4.1:

$\tilde c^*(X)=1-(p_1+p_2)+4b+$

$\phantom{\tilde c^*(X)}+(d_1+d_2)-3b(p_1+p_2)+6b^2+$

$\phantom{\tilde c^*(X)}+0+2b(d_1+d_2)-3b^2(p_1+p_2)+4b^3$.

\noindent We have the following relations:

1. $b\cdot p_i=d_i$ in $IH_2(X)$;

2. $b\cdot d_i=[pt]$ in $IH_0(X)$;

3. $p_i\cdot b_i=0$ and $p_i\cdot b_j=[pt]$ in $H_0(X)$ for $i\not=j$;

\noindent hence

4. the image of $b$ in $IH_4(X)$ is $p_1+p_2$ (from 2 and 3);

5. the image of $b^2$ in $IH_2(X)$ is $d_1+d_2$ (from 1 and 4);

6. $b^2\cdot p_i=[pt]$ in $IH_0(X)$ (from 1 and 2);

7. $b^3=b\cdot (d_1+d_2)=2[pt]$ in $IH_0(X)$ (from 2 and 5).

\noindent We obtain:
$$\tilde c^*(X)=[X]+(3p_1+3p_2)+(4d_1+4d_2)+6[pt]\in IH_*(X)\,.$$
The Chern-Schwartz-MacPherson class is
$$c^*_{MS}(X) = [X]+(3p_1+3p_2)+(4d_1+4d_2)+6[pt]-[vertex]\in IH_*(X)\oplus
IH_*(\{vertex\})\,.$$

\noindent If we compare it with $c^*(X_i)$, then we see that the difference is
supported by $im H_*(\varpi^{-1} \Sigma_X)$ $\subset H_*(\tilde X)$. In
homology, the difference is supported only by the image of $H_*(\Sigma_X)$ in
$H_*(X)$.

Let us come back to the remark we made after Proposition 6.1. 

\proclaim{Observation 7.2.} {The proper inverse image of the cycle $\N 2$ in
the small resolution $X_1$ does not represent the same class as $\tilde
c^2(\tau)\in IH_4(X)$,
so it can not be used to compute the Chern class of $X$.}

\remark{Proof.} We remind that $X_1$ is fibered over $\Pt^1_x$. Since
$\N 2$ is the projective cone over two points in $B$, thus $\tilde \N 2$ is
contained in the
disjoint sum of two fibers in $X_1$. The element $p_2$ is represented
by a fiber. Hence

$[\tilde\N 2]\cdot  p_2=0$.

\noindent We have also

$[\tilde\N 2]\cdot  p_1=2[pt]$,

\noindent see the figures 10--12. This shows, that $[\tilde\N 2]=2
d_2\not=  d_1+ d_2=[\N 2]$.\p\s

If one computes the Chern class of $X$ using $[\tilde\N 2]\in H^4(X_1)$ instead
of $[\N 2]\in IH_2(X_1)$ one obtains the Chern class of $X_1$.
\vfill\eject

\hfil {\bf Cycles in $X$}

The quadric $B$ in $\Pt^3$ with two families of generatrices 1) and 2).

The projective cone $X=cB$ and the generators of
$IH_2(X)$: 3) $d_1$, 4) $d_2$.

The generators of $IH_4(X)$: 5) $p_1$, 6) $p_2$.
\vfill\hfil{\it PICTURES 1 -- 6}
\vfill\eject
\hfil {\bf Cycles in $X$ and $X_1$}

The polar variety $N^ 1$ is the cone over $K$: 7) $N^ 1$, 8) $K$.

The polar variety $N^ 2$ and its proper inverse image in $X_1$: 9)
$N^ 2$, 10) $\tilde N^ 2$.

Proper inverse image in $X_1$ of the generators of $IH_4(X)$:
11) $p_1$, 12) $p_2$.

\vfill\hfil{\it PICTURES 7 -- 12}

\vfill\eject

\baselineskip 14.5pt

\head{References} \parindent 40pt

\item{[BBFGK]} G. \spa{Barthel}, J.-P. \spa{Brasselet}, K.-H. \spa{Fieseler},
O. \spa{Gabber}, L. \spa{Kaup} {\it Rel\`e\-ve\-ment de cycles alg\'ebriques
et homomorphismes  associ\'es en homologie d'intersec\-tion,} Ann. 
Math. 141 (1995), pp. 147-179

\item{[BG]} J.-P. \spa{Brasselet}, G. \spa{Gonzales-Sprinberg} {\it Espaces 
de Thom et contre-exemp\-les de J.L. Verdier et M. Goresky,} Travaux en
cours 23, G\'eom\'etrie Alg\'ebrique et Applications, La Rabida 1984, Hermann
Paris (1987), pp. 5-14

\item{\cite{Ch}} S. S. \spa{Chern} {\it Characteristic classes of hermitian
manifolds} Ann. Math. 47 No. 1 (1946), pp. 85-121

\item{\cite{Eh}} \spa{C. Ehresmann} {\it Sur la topologie de certains
espaces homog\`enes}, Ann. Math. 35 (1934), pp. 396-443

\item{[GM]} M. \spa{Goresky}, R. \spa{MacPherson} {\it Intersection homology
II,} Invent. Math. (1) vol. 72 (1983) pp. 135-116

\item{[Kl]} \spa{Kleiman} {\it Transversality of a general translate,}
Compositio Math. 28 (1974), pp. 287-297

\item{[LT]} \spa{L\^e D\~ung Tr\`ang}, B. \spa{Teissier} {\it Vari\'et\'es
polaires et classes de Chern des vari\'et\'es singuli\`eres} Ann. of Math. 114
(1981), pp. 457-491

\item{[MP1]} R. \spa{MacPherson}~{\it Chern Classes for Singular Algebraic
Varieties}, Annals of Math. 100 (1974), 423 -- 432.\par

\item{[MP2]} R. \spa{MacPherson}~{\it Global questions in the topology of
singular spaces}, Proc. ICM, Warzawa (1983), North Holland 1984, pp. 213 --
235.\par

\item{[MS]} J. \spa{Milnor}, J . \spa{Stasheff} {\it Characteristic Classes,}
Princeton University Press 1974

\item{[Pi]} R. \spa{Piene}, {\it Cycles polaires et classes de Chern pour
les vari\'et\'es projectives singuli\`eres} S\'eminaire Ecole
Polytechnique, Paris, 1977-78 and Travaux en cours 37, Hermann Paris (1988)

\item{[Sc]} M. H. \spa{Schwartz}, {\it Classes caract\'eristiques
d\'efinies par une stratification d'une vari\'et\'e analytique complexe}. C.R.
Acad. Sci. Paris, t. 260, S\'erie I (1965), 3262 -- 3264 et 3535 -- 3537.

\item{[Te]} B. \spa{Teissier}, {\it Vari\'et\'es polaires. II. Multiplicit\'es 
polaires, sections planes et conditions de Whitney.} Algebraic 
geometry (La R\'abida 1981), LNM 961, Springer 1982, 314 --491.

\item{[V]} J. L. \spa{Verdier}, Appendix of [BG].

\item{[We]} A. \spa{Weber}, {\it A morphism of intersection homology
 induced by an algebraic map}. submited.

\item{[Zh]} J. \spa{Zhou}, {\it Classes de Chern pour les vari\'et\'es
singuli\`eres, classes de Chern en th\'eorie bivariante}. Thesis, Marseille,
F\'evrier 1995.

 \end